\documentclass[twoside,11pt]{article}
\usepackage{natbib}
\usepackage{blindtext}

\usepackage{jmlr2e}
\usepackage{amsmath}

\usepackage{lastpage}
\jmlrheading{23}{2022}{1-\pageref{LastPage}}{1/21; Revised 5/22}{9/22}{21-0000}{Umberto Michelucci}

\ShortHeadings{A Singularity Criterion for Countable Gaussian Mixtures}{Umberto Michelucci}
\firstpageno{1}

\begin{document}

\title{A Singularity Criterion for Countable Gaussian Mixtures Based on the Feldman--Hájek Theorem}

\author{\name Umberto Michelucci \email umberto.michelucci@hslu.ch \\
       \addr Computer Science Department\\
       Lucerne Universify of Applied Sciences and Arts\\
       6304, Risch-Rotkreuz, Switzerland}

\editor{My editor}

\maketitle

\begin{abstract}
We study the mutual singularity of countable Gaussian mixture models (GMMs), with particular emphasis on infinite-dimensional settings. We first establish that a countable mixture of Gaussian probability measures is itself a well-defined probability measure. We then prove a general measure-theoretic result showing that if every component of one countable mixture is mutually singular with every component of another, then the two mixtures are mutually singular. Combining this result with the Feldman--Hájek characterization of equivalence and singularity for Gaussian measures yields a sufficient condition for the mutual singularity of countable Gaussian mixtures. We also discuss the mixed case, in which the presence of equivalent components prevents mutual singularity and leads naturally to a decomposition into singular and absolutely continuous parts.
To illustrate these theoretical results, we present a series of numerical experiments involving high-dimensional Gaussian mixture models. The experiments demonstrate the emergence of increasing separability with dimension under different mechanisms, including mean shifts, covariance differences, and independently generated random mixtures. A complementary experiment with a shared Gaussian component shows that complete asymptotic separation fails when the pairwise singularity condition is violated. Together, the theoretical and numerical results provide a measure-theoretic framework for understanding asymptotic separability in high-dimensional Gaussian mixture models.
\end{abstract}

\begin{keywords}
  Gaussian Mixture Models, classification, probability, Infinite-Dimensional Measures
\end{keywords}

\section{Introduction}
The study of probability measures in high-dimensional spaces is often dominated by the phenomenon of measure concentration \citep{billingsley1999, bogachev1998, feldman1958, hajek1958}.  As one transitions to infinite-dimensional settings, the geometry of these measures changes radically, leading to the behavior described for Gaussian measures in 1958 by the Feldman-H\'{a}jek dichotomy \citep{feldman1958, hajek1958}. This result belongs to a broader class of dichotomy theorems for infinite-dimensional measures, the most notable being Kakutani's theorem on infinite product measures \citep{kakutani1948}. The Feldman-Haj\'ek theorem states that two Gaussian measures in an infinite-dimensional Hilbert space are either equivalent or mutually singular.
While the behavior of individual Gaussian measures is well-documented, modern statistical learning frequently employs Gaussian Mixture Models (GMMs) to represent complex distributions, and thus it is a relevant questions under which conditions two different GMMs can be perfectly separable in infinite dimensions. Some general discussion about sums of random variables and mixtures can be found in the book by Vakhania \textit{et al.} in \citep{vakhania1987probability}, but not a specific result for GMMs.
Several authors have studied GMM models in high dimensional settings, but to the best of the author's knowledge,
the literature does not contain an explicit
singularity criterion for countable Gaussian
mixture models in infinite-dimensional settings. Almost all research work found dealt with the problem of classifying or learning GMM models, and were not interested in finding conditions for mutual singularity. 
It may be of interest to the reader to highlight some of the works related to learning GMMs. It should be emphasised that these contributions are in a certain sense connected to and conceptually aligned with the research presented in this work. In fact, if two GMMs are mutually singular, it also means that they are perfectly separable, or in other words, under the assumption that one finds the right classifier, one can classify them with zero missclassification error. A first interesting work in this direction, is the one by Mignacco \textit{et al.} \citep{mignacco_role_2020}. In their work, they studied the role of regularisation in a classification problem of high-dimensional two-mixture of noisy gaussians. Refinetti \textit{et al.} studied the problem of classifying high-dimensional Gaussian mixtures, with a focus on neural networks \citep{refinetti_classifying_2021} and showed the superiority of neural networks over kernel methods. Azizyan \textit{et al.} in \citep{azizyan_minimax_2013} also studied the problem of learning Gaussian mixtures with equal mixing weights in a high-dimensional settings under small-mean separation. Other works goes along the same line, in studying how to \textit{learn} gaussian mixtures. The objective of this study is not to address that point directly, but rather to investigate how the Feldman--Hájek characterization of Gaussian singularity can be combined with a general measure-theoretic result on countable mixtures to obtain a singularity criterion for Gaussian mixture measures. As already mentioned, this work is important for the problem of learning mixtures, as it deals with the limit where error-free (perfect) classification becomes theoretically possible. 
The primary contribution of this work is the derivation of a sufficient condition for mutual singularity of countable Gaussian Mixture Models. More specifically, we establish a general measure-theoretic result showing that countable mixtures inherit mutual singularity from pairwise singular components, and we combine this result with the Feldman--Hájek characterization of Gaussian measures to obtain a singularity criterion for Gaussian mixtures.
Specifically, we establish the following.
(i) We provide a rigorous proof that a countable mixture of Gaussian measures is a well-defined probability measure. (ii) We prove that the mutual singularity of all  Gaussian components of the mixtures is a sufficient condition for the mutual singularity of the resulting mixtures. (iii) We discuss briefly the ``mixed case'' demonstrating how the existence of equivalent components breaks the dichotomy and leads to partial absolute continuity.
By extending these results to GMMs, we provide the first steps for the theoretical foundation for the study of asymptotic separability in complex, multi-modal distributions. (iv) We show, through numerical simulations, how the theoretical results proven in this paper appear in real cases.
The rest of this paper is organised as follows. In Section \ref{sec:formalism}, we establish the formal measure-theoretic framework for countable mixtures and prove that they are valid probability measures via Lemmas \ref{lemma:mu-finite} and \ref{lemma:mu-measure}. Section 3 presents our main result,
a singularity criterion for countable Gaussian mixtures. (Theorem \ref{th:main}), and its proof. Section \ref{sec:exp} describes the experiments, and finally, Section \ref{sec:discussion} discusses the implications of these results.

\section{Formalism and Gaussian Mixtures}
\label{sec:formalism}

Let $H$ be a real separable Hilbert space and let $\mathcal{B}(H)$ denote its
Borel $\sigma$-algebra. Throughout this paper, Gaussian measures are understood
in the sense of Gaussian probability measures on $(H,\mathcal{B}(H))$.

A Gaussian measure $\gamma_k$ on $H$ is denoted by
\[
\gamma_k=\mathcal N(m_k,\Sigma_k),
\]
where $m_k\in H$ is its mean and $\Sigma_k:H\rightarrow H$ is its covariance
operator. The covariance operator is the unique positive, self-adjoint,
trace-class operator satisfying
\[
\langle \Sigma_k u,v\rangle
=
\int_H
\langle x-m_k,u\rangle
\langle x-m_k,v\rangle
\,d\gamma_k(x),
\qquad u,v\in H.
\]

This is the standard characterization of Gaussian measures in
infinite-dimensional Hilbert spaces (see, e.g., \citep{bogachev1998,
vakhania1987probability}).

A countable Gaussian mixture is a probability measure $\mu$ defined for every
$A\in\mathcal B(H)$ by
\begin{equation}
\label{eq:mu-def}
\mu(A)=\sum_{k=1}^{\infty}\pi_k\gamma_k(A),
\end{equation}
where $\gamma_k=\mathcal N(m_k,\Sigma_k)$ are Gaussian probability measures on
$(H,\mathcal B(H))$ and the mixing weights satisfy
\begin{equation}
\pi_k\ge 0,
\qquad
\sum_{k=1}^{\infty}\pi_k=1.
\end{equation}

In Lemma \ref{lemma:mu-measure} we prove that $\mu$ itself is a probability measure. Lemma \ref{lemma:mu-finite} is a stepping stone to prove lemma \ref{lemma:mu-measure}. Note that although the lemma is easily proven, it serves a pedagogical purpose because it isolates the property needed later in the proof, and thus it is reported here in its entirety.
\begin{lemma}
\label{lemma:mu-finite}
Given a Gaussian measure $\gamma_i$, for any sequence of disjoint sets $\{A_n\}_{n=1}^\infty$ with $A_n\in \mathcal{B}(H)$ the series
$
\sum_{n=1}^{\infty}\gamma_i(A_n)
$
converges and its sum is bounded above by $1$.
\end{lemma}
\begin{proof}
    By definition $A_n$ are disjoint sets, and since $\gamma_i$ is a measure then
    \begin{equation}
        \sum_{n=1}^\infty \gamma_i(A_n)=\gamma_i\left( \bigcup_{n=1}^\infty A_n\right)
    \end{equation}
    But also by definition $\cup_{n=1}^\infty A_n \subseteq H$ and therefore
    \begin{equation}
        \gamma_i\left( \bigcup_{n=1}^\infty A_n\right) \leq \gamma_i(X) = 1
    \end{equation}
    Additionally, since each $\gamma_i(A_n)\geq 0$ then the partial sums $s_m=\sum_{n=1}^m \gamma_i(A_n)$ are monotonically increasing and bounded above by 1, hence the partial sums converge. This concludes the proof.
\end{proof}
\begin{lemma}[$\mu$ is a probability measure]
\label{lemma:mu-measure}
    Let us define $\mu$ by 
    \begin{equation}
    \mu(A)=\sum_{k=1}^\infty \pi_k \gamma_k(A)
\end{equation}
where $\gamma_k$ are probability Gaussian measures on $(H, \mathcal{B}(H))$ characterised by the mean vector $m_k$ and by their covariance $\Sigma_k$. $\pi_k$ are mixing weights such that $\pi_k\geq 0$ and $\sum_{k=1}^\infty \pi_k=1$. $\mu$ is a probability measure.
\end{lemma}
\begin{proof}
    To prove that $\mu$ is a measure we need to prove four results from the definition of a probability measure.
    \begin{enumerate}
        \item[(i)] Non-negativity. $\mu(A)\geq 0$ since all $\pi_k\geq0$ and, since they are measures, $\gamma_k(A)\geq 0$.
        \item[(ii)] Empty set. $\mu(\emptyset)=\sum_{k=1}^\infty\pi_k\gamma_k(\emptyset)=\sum_{k=1}^\infty\pi_k \cdot 0=0$.
        \item[(iii)] $\sigma$-additivity. For any sequence of disjoint sets $\{A_n\}_{n=1}^\infty$, we can write
        \begin{equation}
            \mu\left( \bigcup_{n=1}^\infty A_n\right)=\sum_{i=1}^\infty \pi_i\gamma_i \left( \bigcup_{n=1}^\infty A_n\right) =\sum_{i=1}^\infty \pi_i \sum_{n=1}^\infty \gamma_i(A_n)
        \end{equation}
        By Lemma \ref{lemma:mu-finite}, the second sum is finite (the partial sums converge), and since all terms $\pi_i\gamma_i(A_n)$ are nonnegative, Tonelli's theorem permits interchange of the two infinite sums. We can thus exchange the two sums, giving
        \begin{equation}
            \mu\left(\bigcup_{n=1}^\infty A_n\right)=\sum_{n=1}^\infty \sum_{i=1}^\infty \pi_i \gamma_i(A_n) = \sum_{n=1}^\infty \mu(A_n)
        \end{equation}
            
            \item[(iv)]Finally,
\begin{equation}
\mu(X)
=
\sum_{k=1}^{\infty}\pi_k \gamma_k(X)
=
\sum_{k=1}^{\infty}\pi_k
=
1,
\end{equation}
since each \(\gamma_k\) is a probability measure and
\(\sum_{k=1}^{\infty}\pi_k=1\).
Therefore \(\mu\) is a probability measure.
    \end{enumerate}
\end{proof}
Now let us consider two GMMs $\mu$ and $\nu$ defined over two disjoint sets of indices $I_1$ and $I_2$ similarly to Equation \ref{eq:mu-def}. 
Now we will prove the following lemma.
\begin{lemma}
\label{lemma:singularity}
Let $I_1$ and $I_2$ be two disjoint countable index sets.
Let $(\gamma_i)_{i\in I_1}$ and $(\tau_j)_{j\in I_2}$ be Gaussian measures such that
\begin{equation}
\gamma_i \perp \tau_j
\end{equation}
for every $i\in I_1$ and every $j\in I_2$.
Define
\begin{align}
\mu &= \sum_{i\in I_1}\pi_i\gamma_i,\\
\nu &= \sum_{j\in I_2}\rho_j\tau_j,
\end{align}
where
\begin{equation}
\pi_i,\rho_j\ge 0,
\qquad
\sum_{i\in I_1}\pi_i=1,
\qquad
\sum_{j\in I_2}\rho_j=1.
\end{equation}
Then
\begin{equation}
\mu\perp\nu.
\end{equation}
\end{lemma}

\begin{proof}
    To prove this lemma, we need to show that $\exists A$ such that $\mu(A)=1$ and $\nu(A)=0$. Let us construct such a set. Recall that, by assumption, $\gamma_i \perp \tau_j$ (the two measures are mutually singular) $\forall i\in I_1$ and $\forall j \in I_2$, thus 
    there exists a measurable set $A_{ij}$ such that
\begin{equation}
\gamma_i(A_{ij})=1
\qquad\text{and}\qquad
\tau_j(A_{ij})=0.
\end{equation}
    $\forall i\in I_1$ and $\forall j \in I_2$. Consider the set
    \begin{equation}
    \label{eq:A}
        A= \bigcup_{i\in I_1} \left( \bigcap_{j\in I_2} A_{ij}\right)
    \end{equation}
    We want to show that $\mu(A)=1$ and $\nu(A)=0$, thus proving the lemma. To simplify our notation, let us define
    \begin{equation}
        S_i=\bigcap_{j\in I_2} A_{ij}
    \end{equation}
    First of all, let us show that
\begin{equation}
\tau_j(S_i)=0    
\end{equation}
for every $i\in I_1$ and $j\in I_2$.
Since
\begin{equation}
    S_i=\bigcap_{k\in I_2} A_{ik},
\end{equation}
it follows that
\begin{equation}
    S_i \subseteq A_{ij}
\end{equation}
for every fixed $j\in I_2$. Therefore
\begin{equation}
    \tau_j(S_i)
    \leq
    \tau_j(A_{ij})
    =
    0.
\end{equation}    
    Now we can prove that $\gamma_i(S_i)=1$ for $j\in I_2$ and $i\in I_1$.
    To prove this it is enough to prove that $\gamma_i(S_i^c)=0$ since $\gamma_i$ is a probability measure. Thanks to the De Morgan's law
    \begin{equation}
        S_i^c = \left( \bigcap_{j\in I_2} A_{ij}\right)^c = \bigcup_{j\in I_2} A_{ij}^c
    \end{equation}
    it follows
    \begin{equation}
        \gamma_i(S_i^c)\leq \sum_{j\in I_2} \gamma_i (A_{ij}^c) = 0
    \end{equation}
    since the $A_{ij}^c$ are not necessarly mutually disjoint so we use the $\sigma$-subadditivity property of $\gamma_i$. So we established that $\gamma_i(S_i)=1$ for $i\in I_1$ and
$\tau_j(S_i)=0$ for $j\in I_2$ and $i\in I_1$.
    Now let us turn our attention to $A$ defined in Equation (\ref{eq:A}). 
\begin{equation}
    \nu(A)
=
\nu\!\left(
\bigcup_{i\in I_1}S_i
\right)
\le
\sum_{i\in I_1}\nu(S_i)
=
\sum_{i\in I_1}
\sum_{j\in I_2}
\rho_j\tau_j(S_i)
=
0.
\end{equation}
    since we have proved that $\tau_j(S_i)=0$. Thus we have proven that $\nu(A)=0$. 
    Now let us prove that $\mu(A)=1$. Since
\begin{equation}
    A=\bigcup_{i\in I_1}S_i,
\end{equation}
we have $S_i\subseteq A$ for every $i\in I_1$. Since we have already proved that
\begin{equation}
    \gamma_i(S_i)=1,
\end{equation}
it follows that
\begin{equation}
    \gamma_i(A)=1
\end{equation}
for every $i\in I_1$. Therefore
\begin{equation}
    \mu(A)
    =
    \sum_{i\in I_1}\pi_i\gamma_i(A)
    =
    \sum_{i\in I_1}\pi_i
    =
    1.
\end{equation}
Thus we have proven that $\mu(A)=1$.
Thus we have proven that $\mu\perp \nu$.
\end{proof}
We can now state the main singularity criterion for countable Gaussian mixtures.

\section{Mutual Singularity of Countable Gaussian Mixtures}

Based on the construction of the measure $\mu$ as a countable mixture of Gaussian measures the singularity result established in Lemma \ref{lemma:singularity}, we state the main result.

\begin{theorem}[Mutual Singularity of Countable Gaussian Mixtures]
\label{th:main}
Let $H$ be a real separable Hilbert space with Borel
$\sigma$-algebra $\mathcal B(H)$. Let $\mu$ and $\nu$ be two Gaussian Mixture Models (GMMs) defined over disjoint index sets $I_1$ and $I_2$ as
\begin{equation}
    \mu = \sum_{i\in I_1}\pi_i\gamma_i,
\qquad
\nu = \sum_{j\in I_2}\rho_j\tau_j.
\end{equation}
where
\begin{equation}
\gamma_i=\mathcal N(m_i,\Sigma_i),
\qquad
\tau_j=\mathcal N(n_j,\Lambda_j),
    \end{equation}
are Gaussian measures on H, with means
\(m_i,n_j\in H\) and covariance operators
\(\Sigma_i,\Lambda_j\). The measures $\mu$ and $\nu$ are mutually singular if, for every pair
$(i,j)\in I_1\times I_2$, the Gaussian measures $\gamma_i$ and $\tau_j$
satisfy the Feldman--Hájek singularity conditions, and hence
$\gamma_i\perp\tau_j$.

Specifically, $\mu \perp \nu$ if for all $i \in I_1$ and $j \in I_2$, at least one of the following conditions is met:
\begin{enumerate}
    \item If 
    \begin{equation}
                T_{ij}
=
\Sigma_i^{-1/2}
\Lambda_j
\Sigma_i^{-1/2}
-I
    \end{equation} is not a Hilbert-Schmidt operator, i.e., 
    \begin{equation}
        \sum_{n=1}^{\infty} \lambda_n^2 = \infty
    \end{equation}
    where $\lambda_n$ are the eigenvalues of $T$.
    \item The mean difference $(m_i-n_j)$ does not belong to the common
Cameron--Martin space associated with the two Gaussian measures
$\gamma_i$ and $\tau_j$.
\end{enumerate}
\end{theorem}

\begin{proof}
By Lemma 2.3, if $\gamma_i \perp \tau_j$ for all $i \in I_1$ and $j \in I_2$, then $\mu \perp \nu$. Under the standard Feldman-H\'{a}jek theorem, two Gaussian measures are mutually singular if and only if the aforementioned operator or mean conditions are violated. Since these conditions imply $\gamma_i \perp \tau_j$, the singularity of the mixtures follows from the existence of the set $A = \bigcup_{i \in I_1} (\bigcap_{j \in I_2} A_{ij})$ such that $\mu(A)=1$ and $\nu(A)=0$.
\end{proof}
\section{Numerical Simulations}
\label{sec:exp}

To complement the theoretical results, we conducted a series of numerical experiments designed to investigate different mechanisms that can lead to increasing separability between Gaussian mixture models as the ambient dimension grows. In particular, we consider three distinct scenarios: separability induced by differences in component means, separability induced by differences in covariance structure, and separability arising from independently generated random Gaussian mixture models.

Table~\ref{tab:experiments-overview} provides an overview of the experiments considered in this section. The experiments are intentionally synthetic, allowing the source of separation to be controlled and analyzed in isolation. For each experiment, classification accuracy is evaluated as a function of the ambient dimension in order to assess the emergence of high-dimensional separability.

\begin{table}[t]
\centering
\caption{Summary of the numerical experiments used to investigate the effect of increasing dimension on the separability of Gaussian mixture models.}
\label{tab:experiments-overview}
\begin{tabular}{p{2cm} p{4cm} p{7cm}}
\hline
\textbf{Experiment} & \textbf{Name} & \textbf{Description} \\
\hline

1 &
Mean-Induced Separability &
Two Gaussian mixture models with identical covariance matrices and identical mixture weights. Corresponding components differ only by a constant shift in their means.  \\[0.5em]
2 &
Covariance-Induced Separability &
Two Gaussian mixture models with identical component means and identical mixture weights. Corresponding components differ only in their covariance matrices, with covariances $I_d$ and $\sigma^2 I_d$, respectively. \\[0.5em]
3 &
Independent Random Gaussian Mixtures &
Two independently generated Gaussian mixture models with random means and random diagonal covariance matrices. All mixture components are sampled independently from the same distribution, producing distinct mixtures whose separability emerges from high-dimensional geometry. \\[0.5em]
4 &
Shared-Component Gaussian Mixtures &
Two independently generated Gaussian mixture models sharing one identical Gaussian component (same mean and covariance), while all remaining components are generated independently. This experiment illustrates the effect of violating complete pairwise separability. \\

\hline
\end{tabular}
\end{table}

\subsection{Experiment 1: Mean-Induced Separability}

The purpose of this experiment is to investigate how differences in the component means affect the separability of Gaussian mixture models as the ambient dimension increases.
For a fixed dimension $d$, we construct two Gaussian mixture models
\begin{equation}
\mu_d
=
\frac{1}{K}
\sum_{i=1}^{K}
\mathcal N(m_i,I_d)
\end{equation}
and
\begin{equation}
\nu_d
=
\frac{1}{K}
\sum_{i=1}^{K}
\mathcal N(m_i+\delta,I_d),
\end{equation}
where $K=5$ is the number of mixture components and
\begin{equation}
\delta
=
c\,\mathbf{1}_d
=
c(1,\ldots,1),
\end{equation}
with $c>0$ a fixed constant. In the experiments we use
\begin{equation}
c=0.5.
\end{equation}
The component means are generated independently according to
\begin{equation}
m_i \sim \mathcal N(0,I_d),
\qquad
i=1,\ldots,K.
\end{equation}
The covariance matrices are identical for both mixtures and equal to the identity matrix,
\begin{equation}
\Sigma_i
=
\Lambda_i
=
I_d.
\end{equation}
The mixture weights are uniform,
\begin{equation}
\pi_i=\rho_i=\frac{1}{K},
\end{equation}
so that each component contributes equally to the mixture.
For each component, samples from the first mixture are generated according to
\begin{equation}
X \sim \mathcal N(m_i,I_d),
\end{equation}
while samples from the second mixture are generated according to
\begin{equation}
Y \sim \mathcal N(m_i+\delta,I_d).
\end{equation}
Samples generated from $\mu_d$ are assigned class label $1$, while samples generated from $\nu_d$ are assigned class label $2$.
The resulting dataset is
\begin{equation}
\mathcal D
=
\{(x_i,y_i)\}_{i=1}^{2N},
\end{equation}
with
\begin{equation}
y_i \in \{1,2\}.
\end{equation}
The dataset is divided into training and test subsets using a stratified train--test split. A Gaussian Naive Bayes classifier is trained on the training set and evaluated on the test set.
The experiment is repeated for increasing dimensions
\begin{equation}
d
\in
\{10,50,100,200,500,1000\}.
\end{equation}

\subsection{Experiment 2: Covariance-Induced Separability}

The purpose of this experiment is to investigate whether differences in covariance structure alone are sufficient to induce asymptotic separability as the ambient dimension increases.
For a fixed dimension $d$, we construct two Gaussian mixture models

\begin{equation}
\mu_d
=
\frac{1}{K}
\sum_{i=1}^{K}
\mathcal N(m_i,I_d)
\end{equation}

and

\begin{equation}
\nu_d
=
\frac{1}{K}
\sum_{i=1}^{K}
\mathcal N(m_i,\sigma^2 I_d),
\end{equation}

where $K=5$ is the number of mixture components and $\sigma>1$ is fixed. In our experiments we use
\begin{equation}
\sigma = 1.5.
\end{equation}
The component means are shared between the two mixtures and are generated independently according to
\begin{equation}
m_i \sim \mathcal N(0,I_d),
\qquad
i=1,\ldots,K.
\end{equation}
Thus, the two mixtures differ only through their covariance matrices. The mixture weights are uniform,
\begin{equation}
\pi_i=\rho_i=\frac{1}{K},
\end{equation}
so that each component contributes equally to the mixture.
For each component, samples from the first mixture are generated according to
\begin{equation}
X \sim \mathcal N(m_i,I_d),
\end{equation}
while samples from the second mixture are generated according to
\begin{equation}
Y \sim \mathcal N(m_i,\sigma^2 I_d).
\end{equation}
Samples generated from $\mu_d$ are assigned class label $1$, while samples generated from $\nu_d$ are assigned class label $2$.
The resulting dataset is therefore
\begin{equation}
\mathcal D
=
\{(x_i,y_i)\}_{i=1}^{2N},
\end{equation}
with
\begin{equation}
y_i \in \{1,2\}.
\end{equation}
The dataset is divided into training and test subsets using a stratified train--test split. A Gaussian Naive Bayes classifier is then trained on the training set and evaluated on the test set.
The experiment is repeated for increasing dimensions
\begin{equation}
d
\in
\{10,50,100,200,500,1000\}.
\end{equation}

\subsection{Experiment 3: Randomly Generated GMMs}

The purpose of this experiment is to investigate whether two independently generated Gaussian mixture models become increasingly distinguishable as the ambient dimension grows.
For a fixed dimension $d$, we construct two Gaussian mixture models

\begin{equation}
\mu_d
=
\frac{1}{K}
\sum_{i=1}^{K}
\mathcal N(m_i,\Sigma_i)
\end{equation}

and

\begin{equation}
\nu_d
=
\frac{1}{K}
\sum_{j=1}^{K}
\mathcal N(n_j,\Lambda_j),
\end{equation}

where $K=5$ is the number of mixture components.

The mixture weights are uniform,

\begin{equation}
\pi_i=\rho_j=\frac{1}{K},
\end{equation}
so that each component contributes equally to the mixture.
The component means are generated independently according to
\begin{equation}
m_i \sim \mathcal N(0,I_d),
\qquad
n_j \sim \mathcal N(0,I_d),
\end{equation}
for $i,j=1,\ldots,K$.
The covariance matrices are chosen independently as diagonal matrices,
\begin{equation}
\Sigma_i
=
\operatorname{diag}
(u_{i1},\ldots,u_{id}),
\end{equation}
and
\begin{equation}
\Lambda_j
=
\operatorname{diag}
(v_{j1},\ldots,v_{jd}),
\end{equation}
where
\begin{equation}
u_{ik},v_{jk}
\sim
\operatorname{Uniform}(0.5,1.5)
\end{equation}
independently for all indices.
This construction guarantees that the eigenvalues of all covariance matrices remain bounded away from both zero and infinity, independently of the dimension.
For each mixture, $N$ samples are generated according to the standard mixture sampling procedure. First, a component index is selected uniformly from $\{1,\ldots,K\}$, and then a sample is drawn from the corresponding Gaussian distribution.
Samples generated from $\mu_d$ are assigned class label $1$, while samples generated from $\nu_d$ are assigned class label $2$.
The resulting dataset is therefore

\begin{equation}
\mathcal D
=
\{(x_i,y_i)\}_{i=1}^{2N},
\end{equation}

with

\begin{equation}
y_i
\in
\{1,2\}.
\end{equation}

The dataset is randomly divided into training and test subsets using a stratified train--test split.
A Gaussian Naive Bayes classifier is trained on the training set and evaluated on the test set. The classification accuracy is recorded as a function of the ambient dimension $d$.
The experiment is repeated for increasing dimensions

\begin{equation}
d
\in
\{10,50,100,200,500,1000\}.
\end{equation}

\subsection{Experiment 4: Shared-Component Gaussian Mixtures}

The fourth experiment is a modification of Experiment~3. Two independent Gaussian mixture models are generated in the same manner as described above, with the exception that one component is shared exactly between the two mixtures.
More precisely, if
\begin{equation}
\mu_d
=
\frac{1}{K}
\sum_{i=1}^{K}
\gamma_i
\end{equation}

and

\begin{equation}
\nu_d
=
\frac{1}{K}
\sum_{j=1}^{K}
\tau_j,
\end{equation}

then the first component is chosen such that

\begin{equation}
\gamma_1 = \tau_1,
\end{equation}

that is, the two Gaussian measures have exactly the same mean vector and covariance matrix. All remaining components are generated independently according to the procedure described in Experiment~3.

This construction introduces a non-zero overlap between the two mixtures that persists regardless of the ambient dimension. Consequently, unlike Experiments~1--3, perfect asymptotic classification is not expected. The experiment therefore serves as a numerical illustration of the importance of the pairwise singularity condition appearing in Theorem~3.1.

\section{Results and Discussion}
\label{sec:discussion}

The extension of the Feldman-H\'{a}jek theorem to countable Gaussian mixtures  establishes a rigorous framework for understanding the separability of complex distributions in high dimensions. A critical nuance of this result is its application across both finite and infinite-dimensional Hilbert spaces.

In finite dimensions (e.g. $\mathbb{R}^d$), two non-degenerate Gaussian measures on $\mathbb{R}^d$ are always equivalent,
since both are absolutely continuous with respect to Lebesgue measure with
strictly positive densities. However, singularity ($\gamma_i \perp \gamma_j$) can occur in finite dimensions if the measures are degenerate, for instance, or if they are supported on disjoint lower-dimensional affine subspaces. 

The singularity criterion for countable Gaussian mixtures becomes particularly relevant in infinite-dimensional spaces. The equivalence or singularity dichotomy holds for Gaussian measures also in
finite dimensions. What changes in infinite-dimensional settings is that
mutual singularity becomes much more prevalent. This rigorous disjointness of support in infinite dimensions is a characteristic property of Gaussian measures \citep{feldman1958} and is in general valid for product-type measures \citep{kakutani1948} and, as we have demonstrated, scales to Gaussian Mixture Models (GMMs). If the components of the mixtures are singular, a condition almost always satisfied in the infinite-dimensional limit due to shifts in mean or scaling of covariance, the global mixtures $\mu$ and $\nu$ are themselves perfectly separable.

\subsection{Discussion of Experiments 1--3}

\begin{figure}[hbt]
    \centering
    \includegraphics[width=1.0\linewidth]{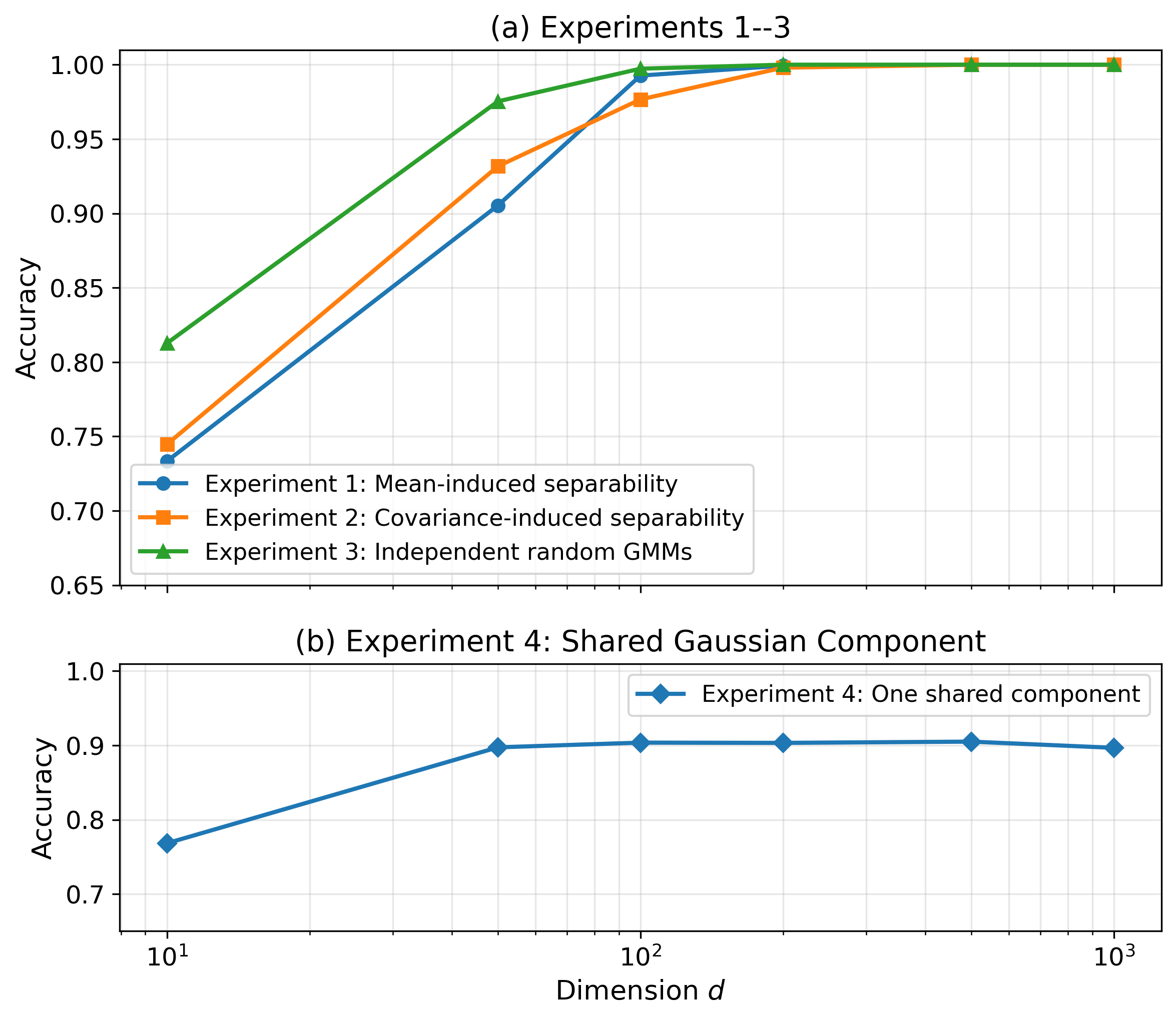}
    \caption{
Classification accuracy as a function of the ambient dimension for the three numerical experiments described in Table~\ref{tab:experiments-overview}. In all cases, the classification accuracy increases with dimension and approaches one, illustrating the emergence of high-dimensional separability between the corresponding Gaussian mixture models.
}
    \label{fig:results1}
\end{figure}

The results of Experiments~1--3 are shown in Figure~\ref{fig:results1}. In all three cases, the classification accuracy increases rapidly with the ambient dimension and approaches one. Although the three experiments are based on different mechanisms (differences in component means, differences in covariance structure, and independently generated random Gaussian mixture models) they exhibit the same qualitative behavior.
These observations are consistent with the main theoretical results of this paper. As the dimension increases, the overlap between the Gaussian mixture models decreases, leading to increasingly accurate classification. In particular, Experiment~3 demonstrates that high-dimensional separability may emerge even when the two mixtures are generated independently from the same underlying distribution of parameters.

It should be emphasized that the experiments are performed in finite-dimensional spaces and therefore do not constitute a direct verification of the mutual singularity results established in Theorem~3.1. Nevertheless, they provide empirical evidence of the same phenomenon: increasing dimension leads to a rapid reduction of overlap between the distributions and consequently to near-perfect classification performance.

\subsection{Discussion of Experiment 4}

The results of Experiment~4 differ  from those of Experiments~1--3. While the classification accuracy increases with dimension, it does not converge to one. Instead, the accuracy approaches approximately $0.9$, as shown in the lower Panel in Figure~\ref{fig:results1}.
This behavior can be explained by the fact that the two Gaussian mixture models share one identical component. Let
\begin{equation}
\mu_d
=
\frac{1}{K}
\sum_{i=1}^{K}
\gamma_i
\end{equation}
and
\begin{equation}
\nu_d
=
\frac{1}{K}
\sum_{j=1}^{K}
\tau_j,
\end{equation}
with
\begin{equation}
\gamma_1 = \tau_1.
\end{equation}
Since the first component is identical in both mixtures, any sample generated from this component has exactly the same distribution under both classes. Consequently, even the optimal Bayes classifier cannot distinguish such samples better than random guessing.

Because the mixture weights are uniform and $K=5$, the probability that a sample originates from the shared component is $1/5$.
Conditional on belonging to the shared component, the probability of classification error is $1/2$.
Therefore, the asymptotic Bayes error is bounded below by
\begin{equation}
\frac{1}{5}\cdot\frac{1}{2}
=
0.1.
\end{equation}
Equivalently, the maximum achievable asymptotic classification accuracy is $1 - 0.1 = 0.9$.
This theoretical prediction is in excellent agreement with the numerical results. As the ambient dimension increases, the remaining four components of the two mixtures become increasingly separated, leading to near-perfect classification on those components. However, the shared component remains intrinsically indistinguishable, producing an irreducible error rate of approximately $10\%$.

From the perspective of Theorem~3.1, this experiment illustrates the necessity of the pairwise singularity condition. The existence of a common Gaussian component implies that the two mixtures possess a non-trivial common part and therefore cannot become completely separated. The observed limiting accuracy of approximately $0.9$ is a direct numerical manifestation of this phenomenon.

\subsection{Implications for Statistical Learning}
The proof that GMMs can be (under certain assumptions) perfectly separable in the limit of infinite dimensions implies the existence of a perfect classifier in the limit. This provides a theoretical ceiling for Bayesian inference and classification tasks involving Gaussian processes and their mixtures, confirming that the complexity of the mixture (even if countable) does not degrade the inherent separability provided by the Feldman-H\'{a}jek conditions.

\subsection{The Mixed Case}

An important aspect in applying the singularity criterion for countable Gaussian mixtures arises when the index sets $I_1$ and $I_2$ each include Gaussian pairs that are not all of the same type, but instead comprise a combination of pairs that are equivalent and pairs that are mutually singular (see again Experiment 4). While the standard theorem for pure Gaussian measures permits only a binary state (equivalence or singularity), the mixture model $\mu$ and $\nu$ allows for a "mixed" or intermediate state.

If there exists at least one pair (i,j)
such that
\begin{equation}
\gamma_i \sim \tau_j ,
    \end{equation}
the perfect separability established in Lemma \ref{lemma:singularity} is lost. In this case, the measures $\mu$ and $\nu$ are no longer mutually singular. Instead, they exhibit partial absolute continuity, which can be formally described via the Lebesgue Decomposition \citep{halmos1974measure}:
\begin{equation}
    \nu = \nu_{ac} + \nu_{s}
\end{equation}
where $\nu_{ac}$ is the portion of the mixture formed by components that share support with $\mu$, and $\nu_{s}$ represents the components that remain singular to all $\gamma_i$ for $i \in I_1$.

In practical terms, this means that the strictly “all-or-nothing” nature of separability in GMMs fundamentally depends on the \textit{pair-wise} singularity condition being satisfied by their individual components. If even a single pair of components in the two mixtures is equivalent, the intersection of their supports creates a region of support overlap for both measures, precluding the existence of an error-free decision boundary. This highlight the importance of the condition $\gamma_i \perp \gamma_j \, \forall i\in I_1,j\in I_2$ in high-dimensional settings.

\section{Practical Applications}

The singularity criterion for countable Gaussian mixtures (Theorem \ref{th:main}) provides a rigorous mathematical foundation for the phenomenon of ``perfect classification'' in high-dimensional settings. While the original theorem is restricted to single Gaussian processes \citep{feldman1958, hajek1958}, real-world datasets, such as high-resolution images, genomic sequences, spectral data, or financial time series, are often characterised by multi-modality and are effectively modeled by Gaussian Mixture Models (GMMs).

In computer science and signal processing, a primary objective is the error-free separation of classes. Theorem \ref{th:main} establishes that if these classes are modelled by GMMs whose components satisfy the Feldman-H\'{a}jek singularity conditions, there exists, under the stated measure-theoretic assumptions, a Bayes classifier
with zero error in the infinite-dimensional limit. This provides a theoretical ceiling for the performance of statistical learning algorithms, confirming that the multi-modal complexity of the mixture does not prevent perfect separability as long as the underlying components are mutually singular. The result shows how pure measure theory revelas itself in applied machine learning by demonstrating that the ``all-or-nothing'' dichotomy remains robust even for countable mixtures.

Since GMMs are universal approximators for continuous probability densities, the implications of Lemma \ref{lemma:singularity} extend to nearly any dataset that can be modelled as a mixture of localised concentrations of measure. In high-dimensional anomaly detection, the theorem predicts that an anomaly, when modelled as a singular mixture component, will not disappear into the baseline data distribution in the infinite-dimensional limit. Instead, the concentration of measure on thin shells ensures that these mixtures remain perfectly distinguishable, as their respective supports are disjoint in the infinite dimensional limit. 

The author  successfully applied this theorem and verified its validity in spectroscopy data (fluorescence spectroscopy of olive oil) and found that in real data, it is possible to reach an extremely high accuracy in classification even from spectral regions that do not contain any chemical information relevant to the classification \citep{michelucci2026} given a sufficient number of input spectral intensities.

For real-world data where complete singularity is not present, the analysis of the ``Mixed Case'' allows researchers to quantify the inherent ``error floor'' of a system. If data analysis reveals that certain components of a GMM are equivalent rather than singular, the Lebesgue decomposition identifies the specific portion of the data that is mathematically inseparable. This allows to distinguish between classification errors caused by sub-optimal model selection and errors that are fundamentally mandatory due to the presence of equivalent probability measures (assuming the number of dimensions is high enough).

\bibliography{references}

\end{document}